\documentclass[preprint,12pt]{elsarticle}

\usepackage{comment}
\usepackage[T1]{fontenc}
\usepackage{lmodern}
\usepackage{amsmath,amssymb,amsthm}
\usepackage{microtype}
\usepackage[hidelinks]{hyperref}
\journal{J. Combin. Theory Ser. B}
\hypersetup{pdftitle={High-chromatic tournament parts without a complete pair of cyclic triangles},
  pdfauthor={Guantao Chen, Hein van der Holst, Shengze Wang}}

\newtheorem{theorem}{Theorem}[section]
\newtheorem{proposition}[theorem]{Proposition}
\newtheorem{lemma}[theorem]{Lemma}
\newtheorem{conjecture}[theorem]{Conjecture}
\DeclareMathOperator{\girth}{girth}
\DeclareMathOperator{\KG}{KG}
\DeclareMathOperator{\conv}{conv}
\newcommand{\B}{\mathcal B}
\newcommand{\R}{\mathbb R}
\begin{document}
\begin{frontmatter}
\title{High-chromatic tournament parts without a complete pair of cyclic triangles}
\author[gsu]{Guantao Chen}
\author[gsu]{Hein van der Holst}
\author[gsu]{Shengze Wang}
\address[gsu]{Department of Mathematics and Statistics, Georgia State University,
Atlanta, GA 30303, USA}

\begin{abstract}
Nguyen, Scott, and Seymour conjectured that any two disjoint subtournaments
of sufficiently large chromatic number contain cyclic triangles, one in
each subtournament, with all edges between them oriented in the same
direction. We disprove this conjecture by constructing, for every positive
integer $d$, a tournament $T$ with a vertex partition $(X,Y)$ such that
$\chi(T[X]),\chi(T[Y])\ge d$, but no cyclic triangle in either part is
complete to a cyclic triangle in the other. Our construction combines
a graph of large girth and large chromatic number with a variant of
Kneser graphs defined in terms of binary functions. We give two proofs
of the required chromatic lower bound, one using the Borsuk--Ulam theorem
and the other using the octahedral form of Tucker's lemma.

\end{abstract}
\end{frontmatter}

\section{Introduction}
\label{sec:intro}

All graphs and tournaments in this paper are finite, and graphs are
simple. For a graph $G$, we write $\chi(G)$ and $\omega(G)$ for its
chromatic number and clique number. For a tournament $T$, its
\emph{chromatic number} $\chi(T)$ is the minimum number of transitive
subtournaments whose vertex sets partition $V(T)$. A tournament is
\emph{transitive} if its vertices admit a linear ordering in which every
arc points forward; equivalently, it contains no directed cycle.
A \emph{cyclic triangle} is a set of three vertices inducing a directed
cycle. A tournament is transitive if and only if it has no cyclic triangle.
For $X\subseteq V(T)$, the notation $T[X]$ denotes the subtournament
induced by $X$. If $X,Y\subseteq V(T)$ are disjoint, we write
$X\Rightarrow Y$ when every arc between them is directed from $X$ to $Y$,
and call $(X,Y)$ a \emph{complete pair}.

The problem of finding complete pairs with large chromatic number is a
central theme of the survey of Nguyen et al.~\cite{NSS}.
Their earlier paper~\cite{NSSElZahar} proposed that, for every integer
$c\ge1$, there exists $d$ such that every tournament $T$ with
$\chi(T)\ge d$ has a complete pair $(A,B)$ satisfying
$\chi(T[A]),\chi(T[B])\ge c$; see also~\cite[Conjecture 1.2]{NSS}.
Klingelhoefer and Newman~\cite{KlingelhoeferNewman2024} showed this is equivalent to a question of El-Zahar and Erd\H{o}s
about anticomplete subgraphs of graphs with bounded clique number. Nguyen et al. provided in their survey paper~\cite{NSS} a short proof, and
obtained a partial result in which $A$ is a
cyclic triangle and $B$ has arbitrarily large chromatic number.

A related line of work asks how the chromatic number of a tournament is
reflected in its out-neighbourhoods. Write
$N_T^+(v)=\{u:v\to u\}$ and $N_T^-(v)=\{u:u\to v\}$.
Berger et al.~\cite{BCCFLSST} conjectured that bounding
$\chi(T[N_T^+(v)])$ for every vertex $v$ bounds $\chi(T)$.
Harutyunyan et al.~\cite{HLTW} proved this in 2019.
Scott and Seymour proposed a strengthening at the 2016 Bellairs
Workshop on Structural Graph Theory: if a graph $G$ and a tournament
$T$ have the same vertex set and $\chi(G)$ is sufficiently large, then
$\chi(G[N_T^+(v)])$ should be large for some $v$; see~\cite[Section 2]{NSS}.
Gir\~ao et al.~\cite{GHILMSS} disproved this strengthening in 2024.
They also constructed tournaments with disjoint sets $X,Y$ of arbitrarily
large chromatic number such that
$$
 \chi(T[Y\cap N_T^+(x)])\le2\quad(x\in X),
 \qquad
 \chi(T[X\cap N_T^+(y)])\le2\quad(y\in Y).
$$
These inequalities still allow cyclic triangles in the cross
out-neighbourhoods. The survey therefore states the following weaker
question, with the explicit observation that its authors suspected it
might be false.

\begin{conjecture}[Nguyen--Scott--Seymour {\cite[Conjecture 2.8]{NSS}}]
\label{conj:NSS}
There exists an integer $d\ge0$ such that, for every tournament $T$ and
all disjoint subsets $X,Y\subseteq V(T)$ with
$\chi(T[X]),\chi(T[Y])\ge d$, there are cyclic triangles
$C_X\subseteq X$ and $C_Y\subseteq Y$ satisfying either
$C_X\Rightarrow C_Y$ or $C_Y\Rightarrow C_X$.
\end{conjecture}

In this paper, we confirm their suspicion by disproving
Conjecture~\ref{conj:NSS}; see Theorem~\ref{thm:main}. The distinction between prescribed parts and the
whole tournament is essential: our result concerns triangles in opposite
prescribed parts, and does not resolve the complete-pair conjecture
of~\cite[Conjecture 1.2]{NSS}.
The construction of Gir\~ao et al.~\cite{GHILMSS} uses Schrijver
subgraphs~\cite{Schrijver} of Kneser graphs. Our construction uses a
variant of Kneser graphs, whose vertices are binary functions
on a graph of large girth and large chromatic number.
Gir\~ao et al. control the out-neighbourhoods of individual vertices
in the other part, while here we control the common in- and
out-neighbourhoods of cyclic triangles.

Section~\ref{sec:construction} defines this graph, states the technical
chromatic lemma, proves the accompanying triangle-free lemma, and
establishes the main construction. In Section~\ref{sec:proofs} we give two proofs of the technical lemma.
The first uses the Borsuk--Ulam theorem, with the geometric terminology
explained explicitly. The second uses only sign vectors and the
octahedral form of Tucker's lemma.

Our main result has been formally verified in Lean~4, assuming
the octahedral Tucker lemma and an Erd\H{o}s graph-existence input.
Details and source code are available at
\url{https://github.com/Swang85/TournamentParts}.

\section{The construction}
\label{sec:construction}

We first state the theorem that disproves Conjecture~\ref{conj:NSS}.

\begin{theorem}
\label{thm:main}
For every positive integer $d$, there is a tournament $T$ with a partition
$V(T)=X\mathbin{\dot\cup}Y$ such that
$$
 \chi(T[X])\ge d\qquad\text{and}\qquad
 \chi(T[Y])\ge d,
$$
but neither $C_X\Rightarrow C_Y$ nor $C_Y\Rightarrow C_X$ holds for
any cyclic triangles $C_X\subseteq X$ and $C_Y\subseteq Y$.
\end{theorem}


Our construction uses a variant of Kneser graphs
associated with a graph $G$, which we call the \emph{binary function
graph} $\B(G)$. Its vertex set is $\{-1,1\}^{V(G)}$. For
$f,g\in V(\B(G))$ and $b\in\{-1,1\}$, let
$$
 I_b(f,g)=\{v\in V(G):f(v)=g(v)=b\}.
$$
Distinct functions $f$ and $g$ are adjacent in $\B(G)$ if and only if
both $I_{-1}(f,g)$ and $I_1(f,g)$ are independent in $G$.
Equivalently, viewing a binary function as the indicator of a subset
with values in $\{-1,1\}$, distinct sets $A,B\subseteq V(G)$ are
adjacent precisely when $A\cap B$ and $V(G)\setminus(A\cup B)$ are
independent in $G$. Thus, instead of requiring the two sets to be disjoint, we require
both their intersection and the complement of their union to be
independent in the base graph \(G\).

In our tournament construction, each function will encode the directions
of the arcs between one vertex and the vertices in the other part.
The two independence conditions will ensure that the common in- and
out-neighbourhoods of vertices corresponding to adjacent functions
are transitive.

The \emph{girth} $\girth(G)$ is the length of a shortest cycle of $G$,
and its \emph{odd girth} is the length of a shortest odd cycle. The
corresponding value is $\infty$ if no such cycle exists. The following
lemma is the technical ingredient of the construction.

\begin{lemma}
\label{lem:topological}
If a graph $G$ has finite odd girth $g$, then $\chi(\B(G))\ge g+1$.
\end{lemma}

We defer the proof to Section~\ref{sec:proofs}, where we give a proof
using the Borsuk--Ulam theorem and a separate combinatorial proof using
Tucker's lemma. Both proofs use the fact that every induced subgraph of
$G$ on fewer than $g$ vertices is bipartite.

\begin{lemma}
\label{lem:trianglefree}
If $\chi(G)>6$, then $\B(G)$ is triangle-free.
\end{lemma}

\begin{proof}
Suppose that $f_1,f_2,f_3$ form a triangle in $\B(G)$. At every vertex
$v\in V(G)$, at least two of the three values $f_1(v),f_2(v),f_3(v)$
agree. Hence the six sets
$$
 I_b(f_i,f_j)\qquad(1\le i<j\le3,\ b\in\{-1,1\})
$$
cover $V(G)$. Every one of these sets is independent in $G$, by the
definition of adjacency in $\B(G)$. Assign each vertex to one of the
sets containing it. The resulting partition is a proper coloring of
$G$ with at most six colors, a contradiction.
\end{proof}


Given a graph $G$ with a linear ordering $<$ of its vertices, its
\emph{backedge tournament} $T(G)$ is defined as follows. For $u<v$,
orient the pair as $v\to u$ if $uv\in E(G)$, and as $u\to v$ otherwise.
Conversely, $G$ is the \emph{backedge graph} of this ordered tournament.
Every independent set of $G$ induces a transitive subtournament of
$T(G)$, because all its arcs point forward. Also, every cyclic triangle
of $T(G)$ contains a backward arc, and therefore contains an edge of $G$.

We use two classical facts, stated in Section~\ref{sec:proofs}:
graphs can have arbitrarily large girth and chromatic number
(Theorem~\ref{thm:erdos}), and if $G$ is a triangle-free backedge graph
of $T$, then $\chi(T)\ge\chi(G)/2$
(Proposition~\ref{prop:backedge}).

\begin{proof}[Proof of Theorem~\ref{thm:main}]
It is enough to prove the result for $d\ge4$, since for $d<4$
we may use the construction for $d=4$.
Fix $d\ge4$. By Theorem~\ref{thm:erdos}, choose a graph $G$ with
$$
 \girth(G)>2d\qquad\text{and}\qquad\chi(G)\ge2d,
$$
and put $B=\B(G)$. Form disjoint sets
$$
 X=\{x_f:f\in V(B)\},\qquad
 Y=\{y_v:v\in V(G)\}.
$$
Choose arbitrary linear orderings on $V(B)$ and $V(G)$.
Give $X$ and $Y$ the corresponding orderings, and let the
subtournaments on these sets be the backedge tournaments of $B$ and
$G$, respectively. For the arcs between the parts, set
$$
 x_f\to y_v\quad\text{if }f(v)=1,
 \qquad
 y_v\to x_f\quad\text{if }f(v)=-1.
$$
These rules specify a tournament $T$ on $X\mathbin{\dot\cup}Y$.

Since $\chi(G)\ge2d>2$, the graph $G$ is not bipartite and has
finite odd girth $g$. Since $\girth(G)>2d$, it follows that $g\ge2d+1$, so
Lemma~\ref{lem:topological} yields
$$
 \chi(B)\ge g+1\ge2d+2.
$$
The graph $G$ is triangle-free, and Lemma~\ref{lem:trianglefree}
shows that $B$ is triangle-free because $\chi(G)\ge8>6$.
Proposition~\ref{prop:backedge} therefore gives
$$
 \chi(T[X])\ge\frac{\chi(B)}2\ge d,
 \qquad
 \chi(T[Y])\ge\frac{\chi(G)}2\ge d.
$$

Let $C_X\subseteq X$ be a cyclic triangle. It contains a backward
arc, so it contains vertices $x_f,x_g$ with $fg\in E(B)$.
If $y_v\in Y$ is dominated by every vertex of $C_X$, then
$f(v)=g(v)=1$. Consequently, the common out-neighbourhood of $C_X$
in $Y$ is contained in $\{y_v:v\in I_1(f,g)\}$.
Similarly, its common in-neighbourhood in $Y$ is contained in
$\{y_v:v\in I_{-1}(f,g)\}$.
Both sets $I_1(f,g)$ and $I_{-1}(f,g)$ are independent in $G$,
since $fg\in E(\B(G))$. Their copies in $Y$ therefore induce
transitive subtournaments. Neither common neighbourhood contains a
cyclic triangle, which excludes both $C_X\Rightarrow C_Y$ and
$C_Y\Rightarrow C_X$ for every cyclic triangle $C_Y\subseteq Y$.
\end{proof}

The proof works for every choice of the two linear orderings.
Moreover, it gives the property that every cyclic triangle in
$X$ has transitive common in- and out-neighbourhoods in $Y$.
Taking $d\ge 4$ disproves Conjecture~\ref{conj:NSS}.

\section{Auxiliary results and two proofs of Lemma~\ref{lem:topological}}
\label{sec:proofs}

We first state the two graph-theoretic results used in the construction.
We then give the two proofs of Lemma~\ref{lem:topological}. The second proof is independent of the geometric
construction in the first.


\begin{theorem}[Erd\H{o}s {\cite{Erdos}}]
\label{thm:erdos}
For all positive integers $a,b$, there exists a graph $G$ such that
$\girth(G)>a$ and $\chi(G)>b$.
\end{theorem}

The following standard inequality is stated and proved
in~\cite[Theorem 4.1]{NSS}. We include its short proof for completeness.

\begin{proposition}
\label{prop:backedge}
If $G$ is the backedge graph of a tournament $T$ under a linear ordering,
then
$$
 \chi(T)\le\chi(G)\le\omega(G)\chi(T).
$$
In particular, if $G$ is triangle-free, then $\chi(T)\ge\chi(G)/2$.
\end{proposition}

\begin{proof}
Every independent set of $G$ is transitive in $T$, giving the first
inequality. For the second, let $W$ induce a transitive subtournament
of $T$. Define a strict partial order on $W$ by $u\prec v$ if $u<v$
in the given linear ordering and $v\to u$ in $T$. Transitivity of
$T[W]$ implies transitivity of this relation. Its comparable pairs
are exactly the edges of $G[W]$, so each chain has at most $\omega(G)$
vertices. Color a vertex by the maximum size of a chain ending at it.
Comparable vertices receive different colors, giving a proper
coloring of $G[W]$ with at most $\omega(G)$ colors. Apply this to
each class in a coloring of $T$, using disjoint palettes.
\end{proof}

\subsection{Proof using the Borsuk--Ulam theorem}
\label{subsec:BU}

For a nonnegative integer $r$, let
$S^r=\{x\in\R^{r+1}:\|x\|=1\}$ be the \emph{unit sphere} in
$\R^{r+1}$. The point $-x$ is the \emph{antipode} of $x$.
We use the following classical form of the Borsuk--Ulam theorem;
see~\cite{Matousek}.

\begin{theorem}[Borsuk--Ulam]
\label{thm:BU}
Let $r,m$ be nonnegative integers with $m\le r$. Every continuous map
$\varphi:S^r\to\R^m$ has a point $x$ such that
$\varphi(x)=\varphi(-x)$.
\end{theorem}

We also recall the elementary geometry of a simplex. The
\emph{convex hull} of points $z_1,\ldots,z_s$ is the set of their
convex combinations $\sum_i a_i z_i$, where $a_i\ge0$ and
$\sum_i a_i=1$. Writing $e_1,\ldots,e_k$ for the standard basis
vectors of $\R^k$, the \emph{standard simplex} is
$$
 \Delta^{k-1}=\conv\{e_1,\ldots,e_k\}
 =\{(a_1,\ldots,a_k):a_i\ge0,\ \textstyle\sum_i a_i=1\}.
$$
Its faces are the convex hulls of subsets of its vertices. In particular,
faces spanned by disjoint sets of vertices are disjoint: a point in both
would have all coordinates zero, contrary to their sum being one.

\begin{proof}[First proof of Lemma~\ref{lem:topological}]
Let $N=|V(G)|$ and $r=g-1$. Since $G$ contains an odd cycle of
length $g$, we have $N\ge g=r+1$. Every induced subgraph on at most
$r$ vertices is bipartite. Let
$c:V(\B(G))\to[k]$ be a proper coloring.
For a binary function $f$, write $(-f)(v)=-f(v)$. The functions
$f$ and $-f$ are distinct and adjacent, since their two common
sign classes are empty. Negation also preserves adjacency.

Choose vectors $a_v\in\R^{r+1}$, indexed by $v\in V(G)$, such that
any $r+1$ of them are linearly independent. Such a choice is explicit:
take distinct real numbers $t_v$ and set
$a_v=(1,t_v,t_v^2,\ldots,t_v^r)$; the required independence is the
Vandermonde determinant identity. The equation $a_v\cdot x=0$
defines a \emph{hyperplane through the origin}, that is, a linear
subspace of dimension $r$. No point of $S^r$ lies on more than $r$
of these hyperplanes, since a vector orthogonal to $r+1$ independent
vectors would have to be zero.

The hyperplanes divide the sphere into \emph{cells}. More explicitly,
a cell consists of all points of the sphere with a fixed pattern of
positive, negative, and zero values for the numbers $a_v\cdot x$.
Thus a cell may be an open region, a boundary piece of a region, or
a lower-dimensional intersection of boundary pieces. These cells and
their boundary relations form a \emph{cell decomposition}: they fit
together to cover the sphere without overlapping interiors.
The decomposition is \emph{antipodally symmetric}, meaning that
$x\mapsto-x$ takes each cell $F$ to another cell $-F$ and preserves
all boundary relations.

Let $\sigma_F\in\{-1,0,1\}^{V(G)}$ be the sign vector of $F$, and
put $Z_F=\{v:\sigma_F(v)=0\}$. We have $|Z_F|\le r$, so $G[Z_F]$
is bipartite. Extend $\sigma_F$ to a binary function $f_F$ by retaining
every nonzero coordinate and assigning the zero coordinates a proper
$2$-coloring of $G[Z_F]$ with values $-1,1$.
Make these choices in antipodal pairs, so that $f_{-F}=-f_F$.
This is possible because negating a proper $2$-coloring is again a
proper $2$-coloring, and $Z_{-F}=Z_F$.

Say that $F$ is a face of $D$ if $F\subseteq\overline D$, where
$\overline D$ is $D$ together with its boundary. In that case, every
nonzero coordinate of $\sigma_F$ agrees with the corresponding
coordinate of $\sigma_D$. Thus $f_F$ and $f_D$ agree outside $Z_F$,
and for each $b\in\{-1,1\}$,
$$
 I_b(f_F,-f_D)
 =\{v:f_F(v)=b,\ f_D(v)=-b\}
 \subseteq\{v\in Z_F:f_F(v)=b\}.
$$
The set on the right is independent in $G$. Also, $f_F\ne-f_D$,
since $|Z_F|<N$ and the two functions $f_F,f_D$ agree at every
vertex outside $Z_F$. It follows that $f_F$ and $-f_D$ are adjacent
in $\B(G)$. By symmetry and negation, whenever $F,D$ are comparable
under the face relation, including $F=D$, we have
\begin{equation}
 c(f_F)\ne c(-f_D).
 \label{eq:comparable}
\end{equation}

We now take an antipodally symmetric \emph{barycentric subdivision}
of this cell decomposition of $S^r$. We describe the construction
explicitly. Choose a point $b_F$ inside each cell $F$, taking its unique
point when the cell consists of a single point. These points are the vertices
of the subdivision. Join $b_{F_0},\ldots,b_{F_t}$ to form a simplex
whenever the distinct cells $F_0,\ldots,F_t$ form a chain under the
face relation. The construction subdivides the cells
successively by dimension, joining the chosen point of each cell to
the already subdivided boundary. The resulting simplices cover the
sphere and meet along common faces; this is what is meant by a
\emph{triangulation} of the sphere.

To obtain antipodal symmetry, choose the points in pairs, with
$b_{-F}=-b_F$. This is possible because no cell equals its antipodal
cell: its nonzero signs are reversed by negation. The face relation is
also preserved by negation, so a chain $F_0,\ldots,F_t$ is paired
with the chain $-F_0,\ldots,-F_t$. Thus every simplex in the
subdivision is accompanied by its antipodal simplex, with the
corresponding vertices paired by $x\mapsto-x$.

For precision, one first forms straight simplices from the points
$b_F$ and then projects them radially onto the sphere, sending a
nonzero vector $z$ to $z/\|z\|$. The vertices associated with a
chain lie in a common open hemisphere, so their convex hull does not
contain the origin and this projection is well defined. Linear
independence also holds: along a strict chain, each new point lies
outside the linear span of the preceding cells. Consequently, radial
projection is one-to-one on each straight simplex, so a point of its
spherical image has well-defined simplex coordinates. Linear
interpolation below is performed in the straight simplex before
projection. Since radial projection commutes with negation, the same
coefficients $\alpha_i$ describe antipodal points using the paired
vertices $b_{F_i}$ and $b_{-F_i}$. These observations explain both
why the map below is well defined on shared faces and why its values
at antipodal points can be compared using the same chain of cells.

Map the new vertex corresponding to $F$ to $e_{c(f_F)}$, and extend
by linear interpolation on each simplex. In detail, if a point $x$
has simplex coordinates $\alpha_0,\ldots,\alpha_t$ in the simplex
corresponding to a chain $F_0,\ldots,F_t$, set
$$
 h(x)=\sum_{i=0}^t\alpha_i e_{c(f_{F_i})},
 \qquad \alpha_i\ge0,\quad\sum_{i=0}^t\alpha_i=1.
$$
The definitions agree on common faces, because they use the same
values at shared vertices. Each interpolation is continuous, so they
combine to give a continuous map $h:S^r\to\Delta^{k-1}$.
Antipodal symmetry gives
$$
 h(-x)=\sum_{i=0}^t\alpha_i e_{c(-f_{F_i})}.
$$
Every pair of cells in the chain is comparable. By
\eqref{eq:comparable}, the colors occurring in $h(x)$ and $h(-x)$
are disjoint. These two points therefore lie in disjoint faces of
$\Delta^{k-1}$, and $h(x)\ne h(-x)$ for every $x\in S^r$.

Finally, all points of $\Delta^{k-1}$ lie in the affine hyperplane
$\sum_i a_i=1$. This hyperplane is a translated copy of
$\R^{k-1}$: its first $k-1$ coordinates determine the last one.
Dropping the last coordinate therefore gives a continuous map from
$S^r$ to $\R^{k-1}$ that still takes different values at every pair
of antipodal points. Theorem~\ref{thm:BU} implies $k-1>r$, and
hence $k\ge r+2=g+1$.
\end{proof}

\subsection{A combinatorial proof using Tucker's lemma}
\label{subsec:Tucker}

Tucker's lemma~\cite{Tucker} is a combinatorial counterpart of the
Borsuk--Ulam theorem. We use its octahedral form
\cite{MatousekKneser,Ziegler}, which is stated entirely in terms of
finite sign vectors and admits a combinatorial proof.
For $\sigma,\tau\in\{-1,0,1\}^N$, write $\sigma\preceq\tau$ if
$\tau_i=\sigma_i$ whenever $\sigma_i\ne0$. Thus $\tau$ is obtained
from $\sigma$ by replacing some zero coordinates by signs. Write
$-\sigma$ for coordinatewise negation and $\mathbf0$ for the all-zero
vector.

\begin{lemma}[Octahedral Tucker lemma]
\label{lem:Tucker}
Let $N\ge2$. If
$$
 \lambda:\{-1,0,1\}^N\setminus\{\mathbf0\}
 \longrightarrow\{\pm1,\ldots,\pm(N-1)\}
$$
satisfies $\lambda(-\sigma)=-\lambda(\sigma)$ for every $\sigma$,
then there are $\sigma\preceq\tau$ such that
$\lambda(\sigma)=-\lambda(\tau)$.
\end{lemma}

The conclusion asserts that two comparable sign vectors receive labels
of equal magnitude and opposite sign. Such a pair is called
\emph{complementary}. No topological terminology is needed to apply
this statement.

\begin{proof}[Second proof of Lemma~\ref{lem:topological}]
Write $V(G)=\{v_1,\ldots,v_N\}$. Suppose, for a contradiction, that
$c:V(\B(G))\to[k]$ is a proper coloring with $k\le g$, and put
$q=N-g\ge0$. For a nonzero sign vector $\sigma$, let
$$
 \operatorname{supp}(\sigma)=\{i:\sigma_i\ne0\},\qquad
 s(\sigma)=|\operatorname{supp}(\sigma)|,
$$
and let $\varepsilon(\sigma)$ be its first nonzero coordinate.
We define a labelling $\lambda$ in two cases.

If $s(\sigma)\le q$, set
$$
 \lambda(\sigma)=\varepsilon(\sigma)s(\sigma).
$$
These labels have magnitudes in $\{1,\ldots,q\}$. This case is
absent when $q=0$.

If $s(\sigma)>q$, then $\sigma$ has at most $g-1$ zero coordinates.
The vertices of $G$ indexed by these coordinates induce a bipartite
graph. Complete $\sigma$ to a binary function $f_\sigma$ by retaining
its nonzero signs and properly $2$-coloring that induced graph on
the zero coordinates. Make the choices in opposite pairs so that
$f_{-\sigma}=-f_\sigma$. The functions $f_\sigma$ and
$-f_\sigma$ are adjacent, and hence have different colors. Set
$$
 \lambda(\sigma)=
 \begin{cases}
 q+c(f_\sigma),&c(f_\sigma)<c(-f_\sigma),\\
 -\bigl(q+c(-f_\sigma)\bigr),&c(-f_\sigma)<c(f_\sigma).
 \end{cases}
$$
Since the smaller of two distinct colors in $[k]$ is at most $k-1$,
these labels have magnitudes between $q+1$ and
$q+k-1\le N-1$. In both cases,
$\lambda(-\sigma)=-\lambda(\sigma)$, so the labelling satisfies
Lemma~\ref{lem:Tucker}.

Choose $\sigma\preceq\tau$ with
$\lambda(\sigma)=-\lambda(\tau)$. The two vectors cannot belong to
different cases, since their label magnitudes would lie in disjoint
ranges. Nor can both satisfy $s\le q$: equal label magnitudes would
give $s(\sigma)=s(\tau)$, and comparability would then force
$\sigma=\tau$, whose nonzero label cannot equal its own negative.
Thus both vectors belong to the second case.

The completions $f_\sigma,f_\tau$ agree at every nonzero coordinate
of $\sigma$. For each sign $b$, the set
$$
 I_b(f_\sigma,-f_\tau)
 \subseteq\{v_i:\sigma_i=0,\ f_\sigma(v_i)=b\}
$$
is independent in $G$, by the choice of the completion on the zero
coordinates. The functions $f_\sigma$ and $-f_\tau$ are distinct,
because $\sigma$ has a nonzero coordinate at which
$f_\sigma=f_\tau$. Hence they are adjacent, and negation also shows
that $-f_\sigma$ and $f_\tau$ are adjacent. Consequently,
$$
 c(f_\sigma)\ne c(-f_\tau),\qquad
 c(-f_\sigma)\ne c(f_\tau).
$$
But opposite labels in the second case force equality in one of these
two pairs, according to which label is positive. This contradiction
proves $\chi(\B(G))\ge g+1$.
\end{proof}

\section*{Acknowledgements}

The mathematical ideas in this paper were developed by the authors.
AI tools were used to improve the writing. The Lean formalization
described in the introduction was generated entirely by AI and
checked by Lean.

\setlength{\bibsep}{3pt}


\begin{thebibliography}{99}

\bibitem{BCCFLSST}
E. Berger, K. Choromanski, M. Chudnovsky, J. Fox, M. Loebl,
A. Scott, P. Seymour, S. Thomass\'e,
Tournaments and coloring,
J. Combin. Theory Ser. B 103 (2013) 1--20.
\href{https://doi.org/10.1016/j.jctb.2012.08.003}{\nolinkurl{doi:10.1016/j.jctb.2012.08.003}}.

\bibitem{Erdos}
P. Erd\H{o}s,
Graph theory and probability,
Canad. J. Math. 11 (1959) 34--38.
\href{https://doi.org/10.4153/CJM-1959-003-9}{\nolinkurl{doi:10.4153/CJM-1959-003-9}}.

\bibitem{GHILMSS}
A. Gir\~ao, K. Hendrey, F. Illingworth, F. Lehner, L. Michel,
M. Savery, R. Steiner,
Chromatic number is not tournament-local,
J. Combin. Theory Ser. B 168 (2024) 86--95.
\href{https://doi.org/10.1016/j.jctb.2024.04.005}{\nolinkurl{doi:10.1016/j.jctb.2024.04.005}}.

\bibitem{HLTW}
A. Harutyunyan, T.-N. Le, S. Thomass\'e, H. Wu,
Coloring tournaments: From local to global,
J. Combin. Theory Ser. B 138 (2019) 166--171.
\href{https://doi.org/10.1016/j.jctb.2019.01.005}{\nolinkurl{doi:10.1016/j.jctb.2019.01.005}}.

\bibitem{KlingelhoeferNewman2024}
F. Klingelhoefer, A. Newman,
Bounding the chromatic number of dense digraphs by arc neighborhoods,
Combinatorica 44 (2024) 881--895.
\href{https://doi.org/10.1007/s00493-024-00098-z}{\nolinkurl{doi:10.1007/s00493-024-00098-z}}.

\bibitem{MatousekKneser}
J. Matou\v{s}ek,
A combinatorial proof of Kneser's conjecture,
Combinatorica 24 (2004) 163--170.
\href{https://doi.org/10.1007/s00493-004-0011-1}{\nolinkurl{doi:10.1007/s00493-004-0011-1}}.

\bibitem{Matousek}
J. Matou\v{s}ek,
Using the Borsuk--Ulam Theorem: Lectures on Topological Methods in
Combinatorics and Geometry,
Universitext, Springer, Berlin, 2003.
\href{https://doi.org/10.1007/978-3-540-76649-0}{\nolinkurl{doi:10.1007/978-3-540-76649-0}}.

\bibitem{NSSElZahar}
T. Nguyen, A. Scott, P. Seymour,
On a problem of El-Zahar and Erd\H{o}s,
J. Combin. Theory Ser. B 165 (2024) 211--222.
\href{https://doi.org/10.1016/j.jctb.2023.11.004}{\nolinkurl{doi:10.1016/j.jctb.2023.11.004}}.

\bibitem{NSS}
T. Nguyen, A. Scott, P. Seymour,
Some results and problems on tournament structure,
J. Combin. Theory Ser. B 173 (2025) 146--183.
\href{https://doi.org/10.1016/j.jctb.2025.02.002}{\nolinkurl{doi:10.1016/j.jctb.2025.02.002}}.

\bibitem{Schrijver}
A. Schrijver,
Vertex-critical subgraphs of Kneser-graphs,
Nieuw Arch. Wisk. (3) 26 (1978) 454--461.

\bibitem{Tucker}
A. W. Tucker,
Some topological properties of disk and sphere,
Proc. First Canadian Math. Congress, Montreal, 1945,
University of Toronto Press, Toronto, 1946, pp. 285--309.

\bibitem{Ziegler}
G. M. Ziegler,
Generalized Kneser coloring theorems with combinatorial proofs,
Invent. Math. 147 (2002) 671--691.
\href{https://doi.org/10.1007/s002220100188}{\nolinkurl{doi:10.1007/s002220100188}}.

\end{thebibliography}
\end{document}